# Graph of NG groups not subset of $S_3$ and $S_4$

Faraj.A.Abdunabi*

Faraj.a.abdunabi@uoa.edu.ly
Frarajarkeas@gmail.com
University of Ajdabyia
*lecturer of Mathematics department

## Abstract

In this paper, we introduce and study the graph of the groups of mappings on a set $X$ with respect to function compositions that cannot be subsets of the symmetric groups $S_3$ and $S_4$. Furthermore, we investigate some of the justifications on these groups that have studied in graph theory, and we got good results. Moreover, we consider the graph of these groups by fixed points and study some properties of it.

## 1-Introduction:

The NG group is group of mappings on a set $X$ with respect to function compositions that cannot be subsets of the symmetric groups called NG-groups or NG-transformations. We try to apply the concepts of the graph theory [1] on these groups on a set X with respect to function compositions that cannot be subsets of symmetric groups S3 and S4. Many researchers study the graph of permutation group [2],[3],[4]. Cheryl E. Praeger study the class consists of all finite transitive permutation groups such that each non-trivial normal subgroup has at most two orbits, and at least one such subgroup is intransitive[5]. Peter Keevash and others determine the maximum number of edges in be the graph on $n$ vertices where two vertices of graphs for all s ≥ 1. We study the graph of the groups of mappings on a set $X$ with respect to function compositions that cannot be subsets of the symmetric groups $S_3$ and $S_4$. Furthermore, we investigate some of the justifications on these groups that have studied in graph theory. Moreover, we consider the graph of these groups by fixed points and study some properties of it. We can generalization of our result on $Sn$ and we introduce the connection between fixed points on NG with graph theory.

## 2-    Preliminaries

In this section, we present some concepts of graph theory. A graph G consist of a nonempty set of vertices and edges (arcs). A simple graph is graph with no loops and no multiple edges. A graph with no edges (i.e. E is empty) is empty and if with no vertices is a null graph. A graph with only one vertex is trivial and the edges are adjacent if they share a common end vertex. Two vertices $u$ and $v$ are adjacent if they are connected by an edge, in other words, $(u, v)$ is an edge. The degree of the vertex $v$, written as $\rho(v)$, is the number of edges with $v$ as an end vertex. By convention, we count a loop twice. If the degree of $u$ is 1 we called a pendant vertex. An isolated vertex is a vertex whose degree is 0. A directed graph (digraph) is formed by vertices connected by directed edges or arcs and we denoted by $(NG)_{dig}$. A simple connected graph is Eulerian if and only if degree of its each vertex is even. Given a graph G, a *walk* in G is a finite sequence of edges of the form $v_0, v_1, \ldots v_n$ also denoted by $v_0 \rightarrow v_x \ldots \rightarrow v_n$ in which any two consecutive edges are adjacent or identical. We call $v_0$ the initial vertex and $v_n$ the end vertex of the walk. The number of edges in a walk called the length of walk. A walk is a trail if any edge is traverse at most once. If $v_0, v_1, \ldots, v_n$ are distinct, then the trail is a path. A path or trail is closed if $v_0 = v_n$ and a closed path containing at least one edge is a circuit. Note that any loop or pair of multiple edges is a cycle. A graph is connected if and only if there is a path between each pair of vertices and is a bipartite graph if and only if each circuit of G has even length. The NG- group is a group of mapping on a finite set that cannot be subsets of symmetric group. In the transformation mapping $f$ on $X$ we called an element $x$ a fixed point if $f(x) = x$, otherwise it is a moved point. We will write the set of all fixed points in NG by $NG_{fix}$ .

**Definition**2-1: The minimum degree of the vertices in a graph G is denoted $\delta$ (G) and the maximum degree of vertices in G we write $\Delta$ (G) and if there is an isolated vertex in G, then $\delta(G) = 0$.

**Definition**2-2:A directed graph is quasi-strongly connected if one of the following conditions holds for every pair of vertices u and v:

*(i)* $u = v$ or

*(ii)* there is a directed path in the digraph from $\underline{u}$ to $v$ or

*(iii)* there is a directed path in the digraph from $v$ to $u$ or

*(iv)* there is a vertex $w$ so that there is a directed path from $w$ to $u$ and a directed path from $w$ to $v$.

Note that quasi-strongly connected digraphs are connected but not necessarily strongly connected.

**Definition**2-3:The vertex $v$ of the digraph G is a root if there is a directed path from $v$ to every other vertex of $G$.

**Definition**2-4: The adjacency matrix of a directed graph NG is $ID = (\rho_{ij})$, where $\rho_{ij}$ = number of arcs that come out of vertex $vi$ and go into vertex $v_j$ .

**Definition 2-5:** A graph is circuitless if it does not have any circuit in it.

**Definition:2-6:** A directed graph with at least one directed circuit is said to be cyclic.

### 3- Graph of NG groups not subset of $S_3$.

In this section, we consider the non-empty set $X = \{a, b, c\}$ with respect to function composition that are not subsets of symmetric groups $S_3$. We can see there 27 mapping from $X$ *to* $X$ and there exits some groups some of them not subsets of Trans($X$) is called NG. We can find the only six N**G-** groups of order two2are: $NG_1=\{(a,a,c),(c,c,a)\}$, $NG_2=\{(a,b,a),(b,a,b)\}$, $NG_3=\{(a,b,b),(b,a,a)\}$, $NG_7=\{(a,c,c),(c,a,a)\}$, $NG_8=\{(b,b,c),(c,c,b)\}$, and $NG_9=\{(b,c,b),\ (c,b,c)\}$.

Note that NG groups are direct graphs and we define the graph of NG as *(NG)_{dig.}*

We can see $NG_1=\{(a,a,c),(c,c,a)\}$ and $NG_7=\{(a,c,c),(c,a,a)\}$ have two fixed pints $NG_{fix}=\{a,\ c\}$ and it is loops . An addition,$(NG_1)_{dig} = (NG_7)_{dig} = (3,6)$.

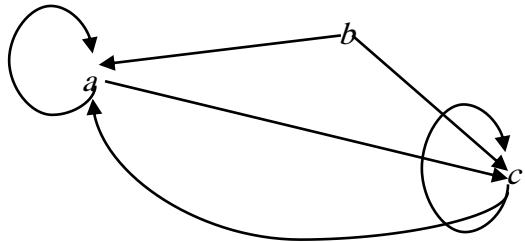

However, $\rho\ (a) = \rho^+\ (a) + \rho^-(a) = 2+3=5$, $\rho\ (b) = \rho^+(b)+ \rho^-(b)=2+0=2$, $\rho\ (c)= \rho^+(c)+ \rho^-(c)=2+3=5$.Therefore, $NG_1$ and $NG_7$are and not symmetric.

By definition 2-1 $\delta$ (NG) =2 and $\Delta$ (NG)= 5.  And By definition 2-4 we can write Adjacency matrix of a $(NG)_{dig}$ as:

$$ID(NG_1)=ID(NG_7)=\begin{pmatrix} \mathbf{1} & \mathbf{0} & \mathbf{1} \\ \mathbf{1} & \mathbf{0} & \mathbf{1} \\ \mathbf{1} & \mathbf{0} & \mathbf{1} \end{pmatrix}$$

Also, the groups $NG_2=\{(a,\ b,\ a),\ (b,\ a,\ b)\}$ and $NG_3=\{(a,b,b),(b,a,a)\}$   are have two fixed points $\{a,\ b\}$ and loops. Moreover, $(NG_2)_{dig} = (NG_3)_{dig} =(3,6)$

However, $\rho\ (a) = \rho^+\ (a) + \rho^-(a) = 2+3=5$, $\rho\ (b)= \rho^+(b)+ \rho^-(b)=2+3=5$, $\rho\ (c)= \rho^+(c)+ \rho^-(c)=2+0=2$.

Adjacency matrix of a $(NG_2)_{dig}$ and $(NG_3)_{dig}$

$$\begin{pmatrix} 1 & 1 & 0 \\ 1 & 1 & 0 \\ 1 & 1 & 0 \end{pmatrix}$$

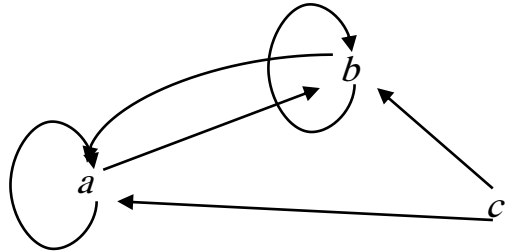

Finally, $NG_8=\{(b,b,c),(c,c,b)\}$ and $NG_9=\{(b,c,b),(c,b,c)\}$ have two fixed pints $\{b, c\}$ also and it is loops. However, $\rho\ (a) = \rho^+\ (a) + \rho^-(a) = 2+0 = 2$, $\rho\ (b) = \rho^+(b) + \rho^-(b) = 2+3 = 5$, $\rho\ (c) = \rho^+(c) + \rho^-(c) = 2+3 = 5$. Also, $(NG_8)_{dig} = (NG_9)_{dig} = (3,6)$

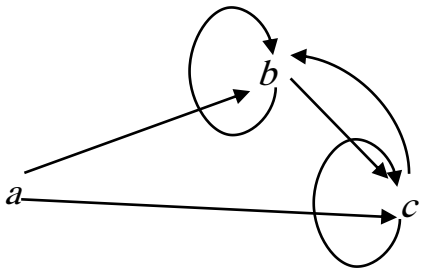

Adjacency matrix of a $(NG_9)_{dig}$ and $(NG_8)_{dig}$

$$\begin{pmatrix} 0 & 1 & 1 \\ 0 & 1 & 1 \\ 0 & 1 & 1 \end{pmatrix}$$

Note that all NG groups of $S_3$ are 1-p directed graphs and $(NG)_{dig} = (3,6)$ and has no isolated vertices and not circuitless since there are loops. Moreover, we can define the adjacency matrix of a directed graph of NG of $S_3$ as:

$$ID(NG)\rho_{ij} = \begin{cases} 0 \text{ if } i \text{ or } j \text{ is movement pint} \\ 1 \qquad\qquad other\ i \neq j \end{cases}.$$

Note that $\rho\ (fix)=5$.

**Handshaking Proposition3-1** :(1736) **[  ]** Suppose G= (*V, E*) graph and the number of vertices is *n* and *m* is the number of the edges. Then $\sum_{v \in V} \rho(v) = 2m$

*Proof*: Since every edge incident with two vertices (different or equals), then every edge contributes exactly two multiply the sum of$(v)$ .Therefore, sum of degree of all vertices is twice the number of edges.

We can apply the handshaking lemma for NG. For example, in NG$_1$ $\sum_{v \in V} \rho(v) = 6 + 2 + 6 = 12 = 2(6)$.

We can rewrite handshaking Proposition3-1as:-

**Proposition 3-2:** Suppose gNG= (*V, E*) the graph of NG groups and the number of vertices is *n* and *m* is the number of the edges. Then$\sum_{v \in V} \rho(v) = 2m$.

**Proposition 3-3:** For every NG has simple circuit between two fixed points.

- It is connected directed graph because it has connected graph but not strongly connected. Since all fixed points, have no directed path to movement points.

By definition2-2 all NG groups of S$_3$ are quasi-strongly connected. And by definition2-3 the movement point in NG S$_3$ is root.

**Proposition 3-**4: Suppose *X* non-empty set of order *n*=3, And G= (*V, E*) graph and number of vertices is n and m the edges then $\sum_{v \in V} \rho(v) = 2n|NG|$.

*Proof:* Let *n*=3. We have $|NG| = 2$ ,So, $\sum_{v \in V} \rho(v) = 2 * 3 * 2 = 12$

**Proposition 3-5**: The graph of all NG of *S$_3$* are bipartite graph and each circuit has even length.

*Proof:* We assume that all NG has no circuit. However, all NG of *S$_3$* has two fixed points. Therefore, NG has circuit, which is contradiction of our assumption. Then the graph of NG is bipartite graph. Suppose NG ={*f,g*}. Since NG is bipartite, we can split it is vertex set into *f* and *g* so that each edge of NG joins a vertex of *f* and a vertex of *g*. Let $v_0 \rightarrow v_1 \rightarrow v_2 \rightarrow v_0$ be a circuit in NG, and assume that $v_0$ is in *f*. Then $v_1$ is in *g*, $v_2$ is in *f*, and so on. Since $v_2$ must be in *g*, the circuit has even length. Moreover, NG of *S$_3$* the only circuit the fixed points and between the fixed points.

**Proposition 3-6:** (Euler Theorem 1736). A connected graph G is Eulerian if and only if the degree of each vertex of G is even.

**Definition: 3-1** : A connected graph is semi-Eulerian if and only if it has exactly two vertices of odd degree.

Note that the graph of NG of $S_3$ is not Eulerian but is it sime- Eulerian.

Note that there is no directed path from any fixed point to the movement points.

Let $S_{fix} = \{v \in V : \exists$ a directed path from any fixed point to other $\}$.

$S_{fix} \neq \emptyset$ as $x \in S_{fix}$ and $S_{fix} \neq V$ as $y \notin S_{fix}$ witch mean $y \in S_{move}$ .

Then $\rho^{-}(S_{move}) = \emptyset$.

**Proposition 3-6**: A digraph $(NG)_{dig}$ of $S_3$ and $S_4$ has at least one root if and only if it is quasi-strongly connected.

*Proof.* If there is a root in the $(NG)_{dig}$ , then it is quasi-strongly connected follows from the definitions 2-2 and 2-3. Suppose $(NG)_{dig}$ ia a quasi-strongly connected and show that it must have at least one root. If $(NG)_{dig}$ is trivial, then it is obvious. Otherwise, we let $V = \{v_1, \ldots, v_n\}$ be the set of of $(NG)_{dig}$ $G$ where $n \geq 2$. The following process shows that there must be a root:

1. Set $P \leftarrow V$ .

2. We take out $v$ from $P$, then we find there exits directed path from distinct vertices $u$ to $v$ in $P$. Equivalently, we set $P \leftarrow P - \{v\}$. We repeat this step as many times as possible.

3. If there is only one vertex left in $P$, then it is the root. For other cases, there are at least two distinct vertices $u$ and $v$ in $P$ and there is no directed path between them in either direction. From condition (*iv*)in definition 2-2 , we have $(NG)_{dig}$ is quasi-strongly connected, it follows that there is a vertex $w$ and a directed path from $w$ to $u$ as well as a directed path from $w$ to $v$. Since $u$ is in $P$, then $w$ can not be in $P$. We take out $u$ and $v$ from $P$ and add $w$, i.e. we set $P \leftarrow P - \{u, v\}$ and $P \leftarrow P \cup \{w\}$. Go back to step #2.

4. Repeat as many times as possible. Every time we do this, there are fewer and fewer vertices in $P$. Eventually, we will get a root because there is a directed path from some vertex in $P$ to every vertex we removed from $P$.

## 4- Graph of NG groups not subset of $S_4$.

We consider $X = \{1, 2, 3, 4\}$ with respect to function composition that are not subsets of symmetric groups, we have the Trans($X$) has 256 mapping from $X$ to $X$.
We have NG-groups of order 2,4,6.

Firstly, the groups of order 2, for example, NG={(1,1,4,4),(4,4,1,1)}. It has two fixed pints {1,4} and it is loops. Moreover, NG is 1-p directed graph and not symmetric with eight edges .i.e $\sum_{v \in V} \rho(v) = 6 + 2 + 2 + 6 = 16 = 2(8)$.

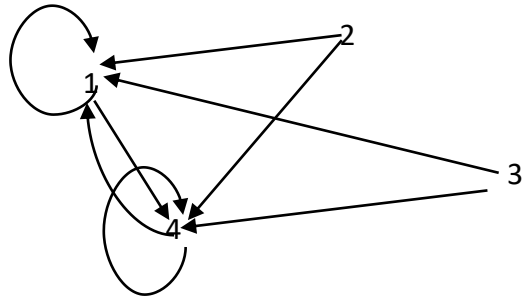

However, $\rho(1) = \rho^+(1) + \rho^-(1) = 2+4=6$, $\rho(2) = \rho^+(2) + \rho^-(2) = 2+0=2$, $\rho(3) = \rho^+(3) + \rho^-(3) = 2+0=2$. $\rho(4) = \rho^+(4) + \rho^-(4) = 2+4=6$. Moreover, $\delta$ (NG)=2 and $\Delta$ (NG)= 6. also, $\rho$ (*fix*)=6.

$$ID(NG) = \begin{pmatrix} 1 & 00 & 1 \\ 1 & 00 & 1 \\ 1 & 00 & 1 \\ 1 & 00 & 1 \end{pmatrix}$$

Secondly, the groups of order 4, for example, NG={(1,1,4,4), (4,4,1,1), (1,4,1,4), (4,1,4,1). It is 2-p directed graph and not symmetric and has two fixed pints {*1,4*} (loops). However, $\rho(1) = \rho^+(1) + \rho^-(1) = 4+8=12$, $\rho(2) = \rho^+(2) + \rho^-(2) = 4+0=4$, $\rho(3) = \rho^+(3) + \rho^-(3) = 4+0=4$. $\rho(4) = \rho^+(4) + \rho^-(4) = 4+8=12$. In addition, $\sum_{v \in V} \rho(v) = 12 + 4 + 4 + 12 = 32 = 2(16)$. $\delta$ (NG)=4 and $\Delta$ (G)= 12.

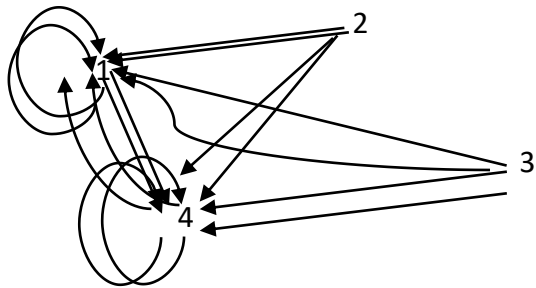

The adjacency matrix of a directed graph NG :

$$ID(NG) = \begin{pmatrix} 2 & 00 & 2 \\ 2 & 00 & 2 \\ 2 & 00 & 2 \\ 2 & 00 & 2 \end{pmatrix} = 2\begin{pmatrix} 1 & 00 & 1 \\ 1 & 00 & 1 \\ 1 & 00 & 1 \\ 1 & 00 & 1 \end{pmatrix}$$

Finally, the groups of order 6, for example, NG= { (1,1,4,3) (1,1,3,4) (4,4,1,3) (4,4,3,1) (3,3,1,4) (3,3,4,1)}. It is 2-p directed graph and not symmetric with three fixed pints $NG_{fix}$={1, 3, 4} are loops. Moreover, there twenty-four edges. $\sum_{v \in V} \rho(v) = 14 + 6 + 14 + 14 = 48 = 2(24)$. And However, $\rho(1) = \rho^+(1) + \rho^-(1) = 6+8=14$, $\rho(2)= \rho^+(2)+ \rho^-(2)=0+6=6$, $\rho(3)= \rho^+(3)+ \rho^-(3)=6+8=14$. $\rho(4)= \rho^+(4)+ \rho^-(4)=6+8=14$. In addition, $\delta$ (NG)=6 and $\Delta$ (G)= 14 and so, $\rho$ (fix)=14.

By definition2-3, we have no root in NG of $S_4$.

We can conclude this proposition

**Proposition 4-**1: Suppose *X* non-empty set of order *n*=4. And G= (*V, E*) graph and number of vertices is *n* and *m* of edges then $\sum_{v \in V} \rho(v) = 2n|NG|$

   *Proof:*

Suppose *n*=4, then we have the three cases of $|NG|$:

Case1:- if $|NG| = 2$, then, $\sum_{v \in V} \rho(v) = 2 * 4 * 2 = 16$ form proposition 3-1.

Case2:- if $|NG| = 4$, then, $\sum_{v \in V} \rho(v) = 2 * 4 * 4 = 32$.

Case3:- if $|NG| = 6$, then, $\sum_{v \in V} \rho(v) = 2 * 4 * 6 = 48$.

                                                                    ∎

**Proposition 4-**3:For every N**G-g**roups of order 2 on a finite set *X, n≥3,* , the difference of number of fixed points of the elements of G is 2. i.e $||f_{fix}| - |g_{fix}|| = 2$

$$\sum_{v \in V} \rho(v) = 2n * ||f fix| - |g fix||$$

*Proof* : By Proposition 4-2 and Proposition3-*1*

**Proposition 4-5**: Every NG-group on $X=\{1, 2, 3\}$, there exits  one of mapping has two fixed points and the image of the third point of this mapping will be one of these fixed points. i.e $f(c)= a\ or\ b$.

*Proof:* Suppose $NG=\{f,g\}$; we have from proposition 3-2 G has six groups . And suppose $f$ has two fixed points. Since $f$ has two fixed points, then we have three cases for $f$ ;

1) if $f(a,b,c)= (a,b,*)$,then $f(c)\neq c$, $f(c)= a\ or\ b$. Then, we have two groups
   2) if $f(a,b,c)= (a,*,c)$, then $f(b)\neq b$, $f(b)= a\ or\ c$. Then we have two groups
   3) if $f(a,b,c)= (*,b,c)$, then $f(a)\neq a$,, $f(a)= b\ or\ c$. Then, we have two groups

From these three cases, we have six NG-groups.

**Proposition 4-6**: Every NG-groups on a finite $X$ has one element having one fixed point.

*Proof:* We assume there exits NG-groups have one element having no fixed points. If we consider the NG-groups of order 2, we have the difference between of number of fixed points of the elements of NG is 2. So, it is contradiction with our assumption.

**Proposition 4-7** : In an acyclic digraph, there exist at least one source (a vertex whose in-degree is zero) and at least one sink (a vertex whose out-degree is zero).

*Proof.* Let G be an acyclic digraph. If G has no arcs, then it is obvious. Otherwise, let us consider the directed path $v_{i0}$ , $e_{j1}$ , $v_{i1}$ , $e_{j2}$ , . . . , $e_{jk}$ , $v_{ik}$ , which has the maximum path length k. Since G is acyclic, $v_{i0} \neq v_{ik}$ . If $(v,\ v_{i0}$ ) is an arc, then one of the following is true:

• $v \neq v_{it}$ for every value of $t = 0$, . . . , k. Then, $v,(v,\ v_{i0}$ ), $v_{i0}$ , $e_{j1}$ , $v_{i1}$ , $e_{j2}$ , . . . , $e_{jk}$ , $v_{ik}$ is a directed path with length k + 1.

• $v = v_{it}$ for some value of $t$. We choose the smallest such $t$. Then, $t > 0$ because there are no loops in G and $v_{i0}$ , $e_{j1}$ , $v_{i1}$ , $e_{j2}$ , . . . , $e_{jt}$ , $v_{it}$ ,$(v,\ v_{i0}$ ), $v_{i0}$ is a directed circuit.

Hence, $\rho^-(v_{i0}$ ) = 0. Using a similar technique, we can show that $\rho^+(v_{ik}$ ) = 0 as well

**Proposition 4-9**: let $X$ be a non-empty set of order $n\geq3$. And G= $(V,\ E)$ graph and number of vertices is n and m the edges then $\sum_{v\in V}\rho(v) \leq 2n!$

We can prove it by *Mathematical induction and proposition 1-1*.

**Proposition 4-10**: Suppose $NG\ of$ order $n$ greaten then 2, then the elements of NG without fixed points is at least $(n-1)*|NG|$.

**Proposition 4-11**: For any NG-group and any fixed point $i \in$ X, $|NG_i| \cdot |i^{NG}| = |NG|$.

**Definition4-2**: Let $(NG)_{dig}$  =(V,E) be dirgraph. For a set $S \subseteq V$. We define $\rho^+(S) = \{w \notin S : \exists v \in$ S s.t.$(v,\ w) \in E)\}$ and $\rho^-(S) = \{w \notin$ S : $\exists v \in$ S s.t.$(w,\ v) \in E)\}$, we can define the NG is strongly connected if there does not exist $S \subseteq V$, $S \neq \emptyset$, $V$ such that $\rho^-(S) = \emptyset$.

**Proposition 4-12:** Every NG is not strongly connected.

*Proof*: If we take S is the set of movement points. Then we can see $\rho^-(S) = \emptyset$.

Then $\rho^-(S_{move}) = \emptyset$. If $z \in \rho^-(S_{move})$ then there exists $w \in \rho^-(S_{move})$ such that $(w, z) \in E$. But then since $w \in \rho^-(S_{move})$ there is a directed path from $x$ to w which can be extended to z, contradicting the fact that $z \notin \rho^-(S_{move})$.

So, NG is not strongly connected.